\documentclass[leqno,10pt,a4paper]{amsart}
\usepackage{amsmath}
\usepackage{amssymb}
\usepackage{latexsym}

\usepackage[all]{xy}
\UseComputerModernTips
\CompileMatrices

\DeclareMathOperator{\Hom}{Hom}
\DeclareMathOperator{\Ext}{Ext}
\DeclareMathOperator{\End}{End}
\renewcommand{\ge}{\geqslant}
\renewcommand{\le}{\leqslant}

\newcommand{\QQ}{\mathcal{Q}}
\newcommand{\OO}{\mathcal{O}}

\newcommand{\N}{\mathbb{N}}

\DeclareMathOperator{\glob}{gl.dim}

\newcommand{\rarr}{\rightarrow}

\DeclareMathOperator{\inj}{i.d.}
\DeclareMathOperator{\proj}{p.d.}

\DeclareMathOperator{\good}{{\mathcal F}(\nabla)}
\DeclareMathOperator{\doog}{{\mathcal F}(\Delta)}

\DeclareMathOperator{\gfd}{\mbox{$\nabla$}.f.d.}
\DeclareMathOperator{\wfd}{\mbox{$\Delta$}.f.d.}

\newcommand{\impl}{\Rightarrow}

\DeclareMathOperator{\rad}{rad}

\newcommand{\Mod}{\mathrm{mod}}

\DeclareMathOperator{\partn}{\Lambda^+}

\newcommand{\GL}{\mathrm{GL}}
\newcommand{\SL}{\mathrm{SL}}
\newcommand{\St}{\mathrm{St}}

\newcommand{\pDelta}{\Delta^\prime}
\newcommand{\pnabla}{\nabla^\prime}

\newcommand{\pP}{P^\prime}
\newcommand{\pT}{T^\prime}
\newcommand{\pS}{S^\prime}
\newcommand{\ple}{\le^\prime}
\newcommand{\op}{^{\mathrm{op}}}

\begin{document}
\theoremstyle{plain}
\numberwithin{subsection}{section}
\newtheorem{thm}{Theorem}[subsection]
\newtheorem{propn}[thm]{Proposition}
\newtheorem{cor}[thm]{Corollary}
\newtheorem{clm}[thm]{Claim}
\newtheorem{lem}[thm]{Lemma}
\newtheorem{conj}[thm]{Conjecture}
\theoremstyle{definition}
\newtheorem{defn}[thm]{Definition}
\newtheorem{rem}[thm]{Remark}
\theoremstyle{remark}
\newtheorem{example}[thm]{Example}

\parskip5pt
\parindent0pt

\title{On the global and $\nabla$-filtration dimensions of 
quasi-hereditary algebras}
\author{Karin Erdmann \and Alison E. Parker}
\address{Mathematical Institute, 24-29 St. Giles',  Oxford OX1 3LB, UK}
\email{erdmann@maths.ox.ac.uk}
\address{School of Mathematics and Statistics,
  University of Sydney, Sydney, NSW 2006, Australia}
\email{alisonp@maths.usyd.edu.au}


\begin{abstract}
In this paper we consider how the $\nabla$--, $\Delta$-- and 
global dimensions of
a quasi-hereditary algebra are interrelated. 
We first 
consider quasi-hereditary algebras with simple preserving duality and
such that 
if $\mu < \lambda$ then $\gfd(L(\mu)) < \gfd(L(\lambda))$ where
$\mu, \lambda$ are in the poset and $L(\mu)$, $L(\lambda)$ are the
corresponding simples. We show that in this case
the global dimension of the algebra is
twice its $\nabla$--filtration dimension.
We then consider more general quasi-hereditary algebras and look at 
how these dimensions are affected by the Ringel 
dual and by two forms of truncation. 
We restrict again to quasi-hereditary
algebras with simple preserving duality and consider various orders on
the poset compatible with quasi-hereditary structure and the $\nabla$-,
$\Delta$- and injective dimensions of the simple and the costandard modules.
\end{abstract}
\maketitle

\section*{Introduction}

Quasi-hereditary algebras were first introduced by Scott \cite{scott} 
in order to
study highest weight categories in the representation theory of semisimple
complex Lie algebras and algebraic groups, and many important results were
proved by Cline, Parshall and Scott (see for example \cite{CPS}). These
algebras can be defined in the context
of arbitrary finite-dimensional algebras, and they were studied from 
this point of view by Dlab and Ringel 
(see for example \cite{dlabring1}, \cite{dlabring2}) and others.
In particular, it turns out that quasi-hereditary algebras are
quite common. 

One important property of quasi-hereditary algebras is that they
have finite global dimension. Furthermore, there is a natural
concept of \emph{$\nabla$-filtration dimension} for representations
of quasi-hereditary algebras. This can be considered as a
generalisation of the notion of injective dimenions. This was introduced for 
algebraic groups by Friedlander and Parshall
\cite{fripar} (where they define the notion of good filtration
dimension which equals our notion of $\nabla$-filtration dimension of a module). 
Later work \cite{parker1,parker2}, 
shows that the understanding
of the $\nabla$-filtration dimension gives a strong hold on
homological properties.

For Schur algebras, the $\nabla$-filtration dimension
of simple modules and of Weyl modules and the projective
dimensions of Weyl modules have nice relationships
with the partial order of the weights. 
Moreover, since there is a duality fixing the simple modules, 
the combination of $\nabla$-filtration dimension and the dual concept
of $\Delta$-filtration dimension, give us an exact relationship between
the $\nabla$-filtration dimension and the injective dimension of a module
(for the regular blocks).

In this paper, we investigate to what extent the interrelations which 
were observed
for Schur algebras 
hold  for arbitrary quasi-hereditary algebras which have a duality fixing
the simple modules. These include the blocks of the category $\OO$
defined by  
{Bern\v ste\u\i n}, Gel'fand and  Gel'fand in  \cite{BGG2}.

 As applications, 
we determine the $\nabla$-filtration dimension and the global dimension of 
the Ringel duals of Schur algebras $S(2,r)$.
Furthermore, we show that a quasi-hereditary algebra
with duality for which the $\nabla$-filtration dimension is strictly
increasing as a function on the poset has
global dimension twice its $\nabla$-filtration dimension.
This applies in particular to the regular blocks of 
Schur algebras $S(n,r)$ with $p>n$ and regular blocks for category
$\OO$ (\cite[theorem 4.7, section 7]{parker2}). 
This proves a particular case of a conjecture of 
Caenepeel and Zhu \cite{CZ} and Mazorchuk and Parker \cite{mazpar}.

\section{Preliminaries}\label{sect:prelim}
\setcounter{subsection}{1}
\begin{defn} \label{defn:qha}
Suppose $S$ is a finite-dimensional algebra over a field $k$.
Let $L(\lambda)$ for $\lambda \in \partn$ be a full set of irreducible
$S$-modules, and let $P(\lambda)$ be the projective cover of
$L(\lambda)$. We  fix a partial order $(\partn, \leq)$.
We then define the \emph{standard module} $\Delta(\lambda)$ to be the largest
quotient of $P(\lambda)$ with composition factors $L(\mu)$ such that 
$\mu \leq \lambda$.

Recall that $S$ is \emph{quasi-hereditary} if
for each $\lambda \in \partn$, 

\begin{enumerate}
\item[(i)]{ the simple module $L(\lambda)$ occurs only once
as a composition factor of $\Delta(\lambda)$, and
}
\item[(ii)]{ the projective $P(\lambda)$ has a filtration by standard modules
where $\Delta(\lambda)$ occurs once, and if $\Delta(\mu)$ occurs then
$\mu \geq \lambda$.
}
\end{enumerate}

The \emph{costandard modules} $\nabla(\lambda)$, are defined dually
by replacing projective by injective modules and quotients by submodules.

We work with finite-dimensional $S$-modules.
We write $\doog$ for the class of $S$-modules
which have a filtration where the sections are $\Delta(\mu)$ for various
$\mu$, and similarly we write $\good$ for the class of $S$-modules
which have a filtration where the sections are $\nabla(\mu)$ for
various $\mu$.
\end{defn}

We henceforth assume that $S$ is a quasihereditary algebra with poset
$(\partn, \leq)$.
Note that by the definition of a quasi-hereditary algebra all the
projective $S$-modules belong to $\doog$ and all the injective
$S$-modules belong to $\good$.
There are other ways of defining quasi-hereditary algebras, but they
turn out to be equivalent. See \cite{klukoe} or
\cite[appendix]{donkbk} for a reasonably self-contained introduction
to quasi-hereditary algebras.

We define $\Ext_S^i(-,-)$ in the usual way (using projective
resolutions) on the category of $S$-modules. We will drop the
subscript if it is clear which category we are working in.

\begin{defn}
Any $S$-module $X$  has a \emph{$\nabla$-resolution}, that
is, there is an exact sequence
$$0 \to X \to M_0 \to M_1 \to \cdots \to M_d \to 0
$$
with $M_i \in \good$.
We say that $X$ has \emph{$\nabla$-filtration dimension} $d$, denoted 
$\gfd(X)=d$ if the following two equivalent conditions hold:

\begin{enumerate}
\item[(i)]{
 $X$ has a $\nabla$-resolution of length $d$ but no $\nabla$-resolution 
of length smaller than $d$;
}
\item[(ii)]{
 $\Ext_S^i(\Delta(\lambda), X) = 0$ for all $i > d$ and
all $\lambda \in \partn$, but there exists
$\lambda \in \partn$ such that $\Ext_S^d(\Delta(\lambda), X) \neq 0$. 
}
\end{enumerate}
(See \cite[proposition 3.4]{fripar} for a proof of the equivalence of
(i) and (ii) where this property is known as the good filtration
dimension of $X$.).
\end{defn}
Dually we have the  notion of \emph{$\Delta$-filtration
 dimension}. This is denoted as $\wfd(X)$.
We also define for a quasi-hereditary algebra $S$,
$$\gfd(S)= \sup \{ \gfd(M) \mid M \mbox{ an $S$-module}\}$$
$$\wfd(S)= \sup \{ \wfd(M) \mid M \mbox{ an $S$-module}\}.$$
But note that $_S S$ considered as a left $S$-module is projective
and hence we have $\wfd(_S S)=0$. Thus we will only use $\wfd(S)$ and
$\gfd(S)$, which are both non-zero in general, as they are defined above.

Recall that $\Ext^i(\Delta (\mu), M)$ for a $S$-module $M$ vanishes
for all $i >0$ and all $\mu \in \partn$ if and only if $M$ has a
$\nabla$-filtration. Thus if $\gfd(M) =0$ then $M \in \good$ and so
the $\nabla$-filtration dimension is a generalisation of this property. 

We also use the notation $\inj(M)$ for the injective dimension of $M$
and $\proj(M)$ for the projective dimension, as well as $\glob(S)$
for the global dimension of $S$.

We have the following important lemma.
\begin{lem}\label{lem:wfd}
{\rm(\cite[lemma 2.2]{parker1}.)} 
For $S$ a quasi-hereditary algebra, $M$, $N$ $S$-modules 
and for $i > \wfd(M) + \gfd(S)$ we have
$$\Ext_S^i(M,N)=0.$$
\end{lem}
As a consequence we have $\glob(S) \le \gfd(S) + \wfd(S)$.

It is possible that different partial orders on the set $\partn$ lead
to the same quasi-hereditary structure. (I.e. different partial orders
may lead to the same standard and costandard modules.)

Once we have a given quasi-hereditary structure (i.e. we are given the
standard and costandard modules) we can replace the
given partial order by a different one 
which gives the same standard
and costandard modules 
but is which more labels would be
incomparable. 

That is if $\lambda < \mu$ and $\lambda$ and $\mu$ are adjacent in the
order 
(that is there is no $\nu \in \partn$ such that
$\mu < \nu < \lambda$), 
but $L(\lambda)$ is not composition factor
of $\Delta(\mu)$ nor of $\nabla(\mu)$ (and hence $\Delta(\lambda)$ is
not a section of $I(\mu)$ by Brauer-Humphreys reciprocity),
then we may safely remove this relation without affecting the
standards or the costandards, since we still get the same modules by
Definition \ref{defn:qha}.

We may continue removing relations in this fashion until we obtain
some minimal partial order which still gives the original standards
and costandards.
Thus, we may assume that if $\mu < \lambda$ and $\mu$ and $\lambda$ are 
adjacent in the order 
then  \label{sect:min}
$L(\mu)$ is a composition factor of $\Delta(\lambda)$ or
of $\nabla(\lambda).$ 

Essentially we have replaced the original parital order by one
that is generated by the preorder
$\mu < \lambda$ if $L(\mu)$ occurs as a composition
factor of $\nabla(\lambda)$ or of $\Delta(\lambda)$.

In this paper we will often assume that $S$ has
a duality $^\circ$ fixing the simple modules. (Such a duality is
sometimes known as \emph{strong duality}.)
For such an algebra, it then follows that
the dual of the costandard module $\nabla(\lambda)^\circ$ is isomorphic
to $\Delta(\lambda)$. It is also clear that $\Ext_S^i(M, N) \cong
\Ext_S^i(N^\circ, M^\circ)$ for all $i\geq 0$, 
and hence that $\gfd(M)=\wfd(M^\circ)$ for
$M, N \in \Mod(S)$.


\section{The global dimension of $S$ with duality}
\label{sect:glob}
Let $S$ be a quasi-hereditary algebra with duality fixing the simple modules.
Then we know that $\glob(S) \leq 2 \gfd(S)$ (as $\wfd(S)=\gfd(S)$
using the remarks above and applying lemma \ref{lem:wfd}). We ask whether
equality holds. (This was orginally conjectured for Schur algebras in
\cite{mythesis} and for more general $S$ in \cite{CZ} and
\cite{mazpar}.)
For most of this section we will be assuming that $S$ satisfies a
particular property which we will call strong property $A$. (We will
weaken this condition slightly in section \ref{sect:props}). That is:
$$\mu < \lambda \impl \gfd(L(\mu)) < \gfd(L(\lambda)).$$
Regular blocks of the Schur algebra satisfy this property
as well as the regular blocks of category $\OO$.
\cite{parker2}.

In the following we write $S \in \mathcal{S}_n$ if
$S$ is a quasi-hereditary algebra with duality fixing 
the simples, with an ordering on the simples such that strong property
$A$ is satisfied and $\gfd(S)=n$.

\subsection{The case with $\nabla$-filtration dimension one}
We first suppose that we have a quasi-hereditary algebra $S$ which belongs
to $\mathcal{S}_1$. 
We split the poset up into a disjoint union
$\partn = \Lambda^+_0\, \dot{\cup}\, \Lambda^+_1$ so that
$\gfd(L(\lambda_0))=0$ for $\lambda_0 \in \Lambda^+_0$ and 
$\gfd(L(\lambda_1))=1$ for $\lambda_1 \in \Lambda^+_1$. 
In this case we know that 
$$\Ext_S^2(L(\lambda_1),L(\lambda_1)) \cong \Hom_S(Q^\circ,Q)\ne 0$$
where $Q=\nabla(\lambda_1)/L(\lambda_1)$ 
using~\cite[lemma 2.6]{parker1}.
So clearly the algebra has $\glob(S)=2$.
For the induction to come we will use the following.

\begin{lem} Let $S$ be in ${\mathcal S}_1$.  
Then for all $S$-modules $Q$ with $\gfd(Q)=1$ we have
$$\Ext_S^2(Q^\circ, Q) \neq 0. $$ 
\end{lem}
\begin{proof} 
We first note that the presence of strong property $A$ 
gives us that if $\lambda_1$ and $\lambda_2$ are both
in $\Lambda^+_i$ for $i \in \{ 0,1\}$ then $\lambda_1$ and $\lambda_2$
are incomparable and hence $\Ext_S^1(L(\lambda_1), L(\lambda_2)) \cong
\Ext_S^1(L(\lambda_2), L(\lambda_1)) \cong 0$. 
This in particular implies that the injective hull of $L(\lambda_i)$ is
$\nabla(\lambda_i)$ for $\lambda_i \in \Lambda^+_1$ and also that
the quotient $\nabla(\lambda_i)/L(\lambda_i)$ is a direct sum of
simples $L(\mu_j)$ with $\mu_j \in \Lambda^+_0$.

Case (1): Assume first that the socle of $Q$ has only
$L(\lambda_i)$ with $\lambda_i \in \Lambda^+_1$.
Then we have a (non-split) injective hull
$$0 \to Q \to I(Q) = \bigoplus_{i} \nabla(\lambda_i) \to N \to 0
$$
and $N$ is a direct sum of copies of $L(\mu_j)$ with $\mu_j \in
\Lambda^+_0$.
Applying $\Hom_S(Q^\circ, -)$ gives
$$\Ext_S^1(Q^\circ, N) \cong \Ext_S^2(Q^\circ, Q)
$$
But $\Ext_S^1(Q^\circ, N) \cong \Ext^1(N^\circ, Q)
= \Ext_S^1(N, Q)$ since $N\cong \bigoplus_j L(\mu_j)$ is self-dual. 
This latter $\Ext$ group 
is non-zero (consider the above exact sequence). 

Case (2): Now suppose $Q$ is arbitrary, then we have an exact
sequence
$$0 \to \bigoplus_{j} L(\mu_j) \to Q \to \bar{Q}\to 0
$$
where $\mu_j \in \Lambda^+_0$, $\bar{Q}\neq 0$ and 
has only $L(\lambda_i)$ with $\lambda_i \in
\Lambda^+_1$ in the socle. 
Now $\gfd(L(\mu_j))=0$ and $\gfd(Q)=1$ hence $\gfd(\bar{Q})=1$,
using \cite[lemma 2.5]{parker1}.
Using case (1) we know that $\Ext_S^2(\bar{Q}^\circ, \bar{Q})\neq 0$. 
We will show that there is an epimorphism from
$$\Ext_S^2(Q^\circ, Q) \to \Ext_S^2(\bar{Q}^\circ, \bar{Q})
$$
and this will be enough to show that the first $\Ext$ group is
non-zero.

Apply $\Hom_S(Q^\circ, -)$ to the exact sequence for $Q$, this gives
an exact sequence
$$\ldots \Ext_S^2(Q^\circ, Q) \to \Ext_S^2(Q^\circ, \bar{Q})
\to 0
$$
as  $\bigoplus_j L(\mu_j)$ has injective dimension $\leq 2$.

Now apply $\Hom_S(-,\bar{Q})$ to the exact sequence
$$0 \to \bar{Q}^\circ\to Q^\circ\to \bigoplus_j L(\mu_j)\to 0
$$
This gives an exact sequence
$$\Ext_S^2(Q^\circ, \bar{Q})\to \Ext_S^2(\bar{Q}^\circ, \bar{Q})
\to 0
$$ 
as $\bigoplus_j L(\mu_j)$ has projective dimension $\leq 2$.
The composite of these two maps gives the desired epimorphism.
\end{proof}

\subsection{} Assume now that $S$ is a quasi-hereditary 
algebra in $\mathcal{S}_n$. 
We know an algebra $S_1$ in $\mathcal{S}_1$ has global dimension $2$. 
Moreover for every $S_1$-module $Q$ with $\gfd(Q)=1$ we know
$\Ext_{S_1}^2(Q^\circ, Q) \neq 0$.

\begin{thm}\label{thm:glob}
An algebra $S$ in $\mathcal{S}_n$ has global dimension
$2n$. Moreover for every $S$-module $Q$ with $\gfd(Q)=n$
we have
$$\Ext_S^{2n}(Q^\circ, Q) \neq 0.$$
\end{thm}
\begin{proof}
We have already proved this for $n=1$. We now assume inductively that
an algebra $S_{n-1} \in \mathcal{S}_{n-1}$ has global dimension
$2(n-1)$ and that for every $S_{n-1}$-module $Q$ with $\gfd(Q)={n-1}$
we have
$$\Ext_{S_{n-1}}^{2(n-1)}(Q^\circ, Q) \neq 0.$$

We first show that $\glob(S)=2n$.
We know that $\glob(S) \leq 2\gfd(S)=2n$.
So it is enough to show that $\Ext_{S}^{2n}(L(\lambda), L(\lambda)) \ne
0$,
for $\lambda$ with $\gfd(L(\lambda))=n$.
We have an exact sequence
$$0 \to L(\lambda) \to \nabla(\lambda) \to Q \to 0
$$
and $\gfd(Q)=n-1$. Let $S_{n-1}$ be the quotient
$S/Se_{\Gamma}S$ of $S$ where $\Gamma=\{ \mu \mid
\gfd(L(\mu))=n \}$, then $S_{n-1}$ belongs to 
$\mathcal{S}_{n-1}$. (See section \ref{sect:trunc1} for more details about
$S/Se_{\Gamma}S$.) Moreover $Q$ is an $S_{n-1}$-module, by the
assumptions on $S$.

Applying $\Hom_S(L(\lambda), -)$ to the above exact sequence and then
$\Hom_S(-,Q)$ 
to its $^\circ$-dual, gives us
$$\Ext_S^{2n}(L(\lambda), L(\lambda))
\cong \Ext_S^{2n-1}(L(\lambda), Q) \cong
\Ext_S^{2(n-1)}(Q^\circ,Q)$$  
since $\nabla(\lambda)$ is injective. 
So by the inductive hypothesis we get that 
$$\Ext_S^{2n}(L(\lambda),
L(\lambda))\cong \Ext_S^{2(n-1)}(Q^\circ, Q)
\cong \Ext_{S_{n-1}}^{2(n-1)}(Q^\circ, Q)\neq 0.$$
(For the last equality see section \ref{sect:trunc1}).

Now let $Q$ be some $S$-module with $\gfd(Q)=n$. We must show that
$\Ext_S^{2n}(Q^\circ, Q)\neq 0$.
We note that the modules $\nabla(\lambda_i)$ with $\gfd(L(\lambda_i))=n$
must be injective.
Also $Q$ must have at least one $L(\lambda)$ as a composition
factor with $\gfd(L(\lambda))=n$.

Assume first that the socle of $Q$ is a direct sum of $L(\lambda_i)$'s
with $\gfd(L(\lambda_i))=n$.
Then we have the injective hull
$$0 \to Q \to I \to R \to 0
$$
where $I \cong \oplus \nabla(\lambda_i)$. 
Moreover $\gfd(R)=n-1$ since $I$ is
injective.
We also know that $R$ is an $S_{n-1}$-module, ($S_{n-1}$ as before)
by construction.
So by the inductive hypothesis we know that 
$\Ext_{S_{n-1}}^{2(n-1)}(R^\circ, R)\neq 0$. 
We now have  
$$\Ext_S^{2n}(Q^\circ, Q) \cong \Ext_S^{2(n-1)}(R^\circ, R)
\cong \Ext_{S_{n-1}}^{2(n-1)}(R^\circ, R)
$$
as before as $I$ is injective.

We now consider the general case. 
Let $U \subset Q$ be the largest submodule with 
no composition factors $L(\lambda_i)$ with $\gfd(L(\lambda_i))=n$, 
and let $V = Q/U$. Then 
the socle of $V$ has only composition factors $L(\lambda_i)$, note also that
$V\neq 0$. Consider the exact sequence
$$0 \to U \to Q \to V \to 0.
\leqno{(*)}$$
We know that $\gfd(U)\leq n-1$ since $U$ is an $S_{n-1}$-module.
But $Q$ has $\gfd(Q)=n$ and since $\gfd(V)\leq n$
it follows that $\gfd(V)=n$.
So we know from the first case that $\Ext_S^{2n}(V^\circ, V) \neq 0$.
Therefore it is enough to show that there is
an epimorphism from $\Ext_S^{2n}(Q^\circ, Q)$ onto
$\Ext_S^{2n}(V^\circ, V)$. 

Apply $\Hom_S(Q^\circ, -)$ to the exact sequence $(*)$, this gives
$$\to \Ext_S^{2n}(Q^\circ, Q) \stackrel{\phi}\to \Ext_S^{2n}(Q^\circ, V)
\to \Ext_S^{2n+1}(Q^\circ, U).
$$
The last term is zero since $S$ is known to have global dimension $2n$.

Next, apply $\Hom_S(-, V)$ to the exact sequence
$$0 \to V^\circ\to Q^\circ\to U^\circ\to 0
$$
which gives
$$\to \Ext_S^{2n}(Q^\circ, V)\stackrel{\psi}\to \Ext_S^{2n}(V^\circ, V) \to
\Ext_S^{2n+1}(U^\circ, V)=0.
$$
The composite $\psi \circ \phi$ gives the required epimorphism.
\end{proof}

\subsection{}
The previous section proves a special case of the conjecture of \cite{CZ, mazpar}
that the global dimension of
any quasi-hereditary algebra $S$ with simple preserving duality is twice its
$\nabla$-filtration dimension.
We suspect that a stronger property may be true. That is that one of
the equivalent conditions of the following lemma hold.

\begin{lem}
$\Ext^{2i}(M^\circ, M) \ne 0$ for all $i \le \gfd(M)$
if and only if $\Ext^2(M^\circ,M) \ne 0$ for all $M$ with $\gfd(M)\ne0$.
\end{lem}
\begin{proof}
($\Rightarrow$) clear.
($\Leftarrow$)
Clearly $\Hom(M^\circ, M) \ne 0$ for all $M\ne 0$ as the head of
$M^\circ$ is isomorphic to the socle of $M$.

Now take an injective resolution for $M$ with $d=\gfd(M) \ne0$ (so $M$
is not injective).

$$0\rarr M \rarr I_0 \rarr I_1 \rarr \cdots \rarr I_{d-1} \rarr I_d
\rarr \cdots $$
We denote the images of the map $I_i \rarr I_{i+1}$ by $N_{i+1}$.
We have $\gfd(N_i) = \sup\{0,\ \gfd(M) - i\mid  i \in \N \}$ 
by dimension shifting.

Now suppose $i \le \gfd(M) =d$.
 By dimension shifting and duality we have that
\begin{eqnarray*}
\Ext^{2i}(M^\circ, M) &\cong \Ext^{i+1}(M^\circ, N_{i-1})\\
&\cong \Ext^{i+1}(N_{i-1}^\circ, M)\\
&\cong \Ext^2(N_{i-1}^\circ, N_{i-1})
\end{eqnarray*}
which is non-zero as $i \le d = \gfd(M)$ and so $\gfd(N_i-1) = d -i +1
\ge 1$.
\end{proof}

In a similar vein we have:
\begin{lem}\label{lem:gfd}
$\Ext^{2d}(M^\circ, M) \ne 0$ for $d = \gfd(M)$
if and only if $\Ext^2(M^\circ,M) \ne 0$ for all $M$ with $\gfd(M)=1$.
\end{lem}
Indeed, we have proved 
that the first condition of this lemma holds for 
our special case in the previous section. 
We will give another example of a quasi-hereditary
algebra for which the first condition of this lemma holds 
in example \ref{eg:propA}.

\section{$\nabla$-filtration and Global dimensions for Ringel duals}
In this section we investigate the relationship between the 
$\nabla$-and $\Delta$-filtration dimensions for a quasi-hereditary algebra
and its Ringel dual (as defined in \cite{ringel}).

\subsection{}
A \emph{tilting module} is a module with both a $\nabla$-filtration and a
$\Delta$-filtration. There is a unique indecomposable tilting module
$T(\lambda)$ for each $\lambda \in \partn$ such that $L(\lambda)$
occurs only once and any other composition factor $L(\mu)$ of $T(\lambda)$ has
$\mu < \lambda$. Every tilting module is a direct sum of $T(\mu_i)$ for
some $\mu_i \in \partn$. A \emph{full} tilting module $T$ is a tilting
module for which for all $\mu \in \partn$,
$T(\mu)$ is a direct summand.  
We take a full tilting module $T$ and form a Ringel dual
$\pS=\End_S(T)\op$.
A Ringel dual is also a quasi-hereditary
algebra with poset $(\partn, \ple)$, where $\ple$ is the opposite
ordering to $\le$ on $\partn$. We distinguish the standards, costandards
etc. for a Ringel dual from that of the starting algebra by a prime. 
Different $T$ lead to different `Ringel
duals' but it is unique up to Morita equivalence. So we often say
\emph{the} Ringel dual.
There is a left exact functor $F:S \rarr \pS$ which takes a module
$M$ to $\Hom_S(T,M)$ regarded as an $\pS$-module in the usual manner.

The following relationships hold between various modules for
$S$ and $\pS$.
$\pDelta(\lambda)=F\nabla(\lambda)$, 
$\pP(\lambda)=FT(\lambda)$ and  
$\pT(\lambda)=FI(\lambda)$ for $\lambda \in \partn$.

\begin{propn}\label{propn:ringel}
We have the following equalities.
\begin{enumerate}
\item[(i)]{
$\wfd(\nabla(\lambda)) = \proj(\pDelta(\lambda))$
}
\item[(ii)]{
$\inj(\nabla(\lambda)) = \gfd(\pDelta(\lambda))$
}
\item[(iii)]{
$\proj(\Delta(\lambda)) = \wfd(\pnabla(\lambda))$
}
\item[(iv)]{
$\gfd(\Delta(\lambda)) = \inj(\pnabla(\lambda))$
}
\end{enumerate}
\end{propn}
\begin{proof}
(i). We take a minimal length \emph{tilting resolution} for $\nabla(\lambda)$
$$0 \rarr T_d \rarr \cdots \rarr T_1  \rarr T_0
 \rarr  \nabla(\lambda) \rarr 0$$
 
using \cite[proposition A4.4]{donkbk}. (That is we have a resolution
of shortest possible length
where each $T_i$ is a tilting module). Such a resolution is also a
$\Delta$-resolution for $\nabla(\lambda)$ and hence $d \ge
\wfd(\nabla(\lambda))$. But if the resolution is minimal then
$\Ext^d(\nabla(\lambda), T_d) \ne 0$, thus $d \le
\wfd(\nabla(\lambda))$ as $T_d \in \good$. So $d =
\wfd(\nabla(\lambda))$. We also note that $\Ext^d(\nabla(\lambda),T_d)
\cong \Ext^d(\pDelta(\lambda), FT_d) \ne 0$ using \cite[proposition
A4.8]{donkbk}. So $\proj(\pDelta(\lambda)) \ge d$. 
We now form an projective resolution for $\pDelta(\lambda)$ using
the fact that $F$ is exact on $\good$ \cite[statement (1)(i)
preceeding lemma A4.6]{donkbk}.
$$0 \rarr \pP_d \rarr \cdots \rarr \pP_1  \rarr \pP_0
 \rarr  \pDelta(\lambda) \rarr 0$$
where the $\pP_i = F T_i$ are projective.
Thus $\proj(\pDelta(\lambda)) =d = \wfd(\nabla(\lambda))$.

(ii). We similarly take a minimal length injective resolution for
$\nabla(\lambda))$ and apply $F$  to get a minimal length tilting resolution
for $\pDelta(\lambda)$. By a similar argument to that above we know 
that the length of a minimal tilting resolution for $\pDelta(\lambda)$
is the same as its $\nabla$-filtration dimension.

(iii) and (iv) follow by applying (i) and (ii) to the modules for
$\pS$ and using the fact that $S$ and  $S^{\prime\prime}$ are 
Morita equivalent as quasi-hereditary algebas. 
\end{proof}

\begin{cor}\label{cor:ringel}
Let $S$ be a quasi-hereditary algebra, then
\begin{enumerate}
\item[(i)]{$\gfd(S)=\wfd(\pS)$ and}
\item[(ii)]{$\wfd(S)=\gfd(\pS)$.}
\end{enumerate}
\end{cor}
\begin{proof}
We prove the first statement, the second is similar.
Since 
$$\gfd(S)=\sup\{\proj(\Delta(\lambda))\mid \lambda \in\partn\}$$
we have 
$$\gfd(S)=\sup\{\wfd(\pnabla(\lambda))\mid \lambda \in\partn\}$$
which equals $\wfd(\pS)$ using \cite[lemma 2.10]{parker2}.
\end{proof}

If $S$ has a simple-preserving duality then so does its Ringel dual,
(using \cite[theorem 1]{mcninch} in the case where the induced
automorphism is the identity map).
Hence in this situation we have
$\gfd(S)=\wfd(S)=\wfd(\pS)=\gfd(\pS)$.

\begin{example}
We can now write down various formulae for the $\nabla$-filtration
dimensions of the simple modules for the regular blocks of 
the Schur algebras and their Ringel duals.

Recall that the simples for the Schur algebra, $S(n,r)$, as defined in
\cite{green}, are indexed by the set
of partitions of $r$ into less than or equal to $n$ parts,
$\partn(n,r)$.
A partition $\lambda=(\lambda_1, \lambda_2,\ldots, \lambda_n)$ 
is \emph{regular} if $\lambda_i -\lambda_j \not \equiv i-j \pmod p$
for all $1\le i<j\le n$. 
We say a block is regular if all the partitions in a block are
regular.
We define for $\lambda \in \partn(n,r)$
$$d(\lambda)= \sum_{i=1}^{n-1}\sum_{j=i+1}^{n}
               \left\lfloor \frac{\lambda_i -\lambda_j
                        -i+j-1}{p}\right\rfloor.$$
Now if $\lambda$ is a partition in a regular block with maximal
element $\mu$ we have
$d(\lambda)=\gfd(\Delta(\lambda))=\inj(\pnabla(\lambda))
$
and $d(\mu)-d(\lambda)=\inj(\nabla(\lambda))=\gfd(\pDelta(\lambda))$
using \cite{parker2} and proposition \ref{propn:ringel}.
%
%
%
%
%
%
%
%
\end{example}

\section{Truncation properties}

In this section we investigate the 
behaviour of $\nabla$-filtration dimensions and injective dimensions under
two forms of truncation. Let $S$ be a quasi-hereditary algebra with
poset $\partn$.
We fix a \emph{saturated} subset $\Pi$ of $\partn$, that is $\Pi$ is a subset
of $\partn$ with the property that if $\lambda \in \Pi$  and $\mu \in
\partn$ then $\mu < \lambda$ implies $\mu \in \Pi$. 
We write $\Gamma$ for $\partn \setminus \Pi$. 
Let $e_{\Gamma} := \sum_{\lambda \in \Gamma} e_{\lambda}$, an idempotent;
where we use a fixed 
decomposition of $1$ into a sum of orthogonal primitive
idempotents, and where $e_{\lambda}S$ is a projective module isomorphic
to $P(\lambda)$.

\subsection{}\label{sect:trunc1}
The algebra $S$ has a quotient $S/Se_{\Gamma}S$, denoted
by  
$S(\Pi)$ in \cite[A3.9]{donkbk}. It is quasi-hereditary with respect to
$(\Pi, \leq)$ with standard modules  $\Delta(\lambda)$ and costandard
modules $\nabla(\lambda)$, the same as 
for $S$ when $\lambda \in \Pi$. 

For $M$, $N$ in $\Mod(S(\Pi))$ considered in the natural way 
as a subcategory of $\Mod(S)$,
we have $\Ext^i_{S(\Pi)}(M,N) \cong \Ext^i_S(M,N)$
(see \cite[A3.3]{donkbk} or  \cite[Appendix]{dlabring1}). 
We have the following:
if $M \in \Mod(S(\Pi))$
and if $\Ext_S^i(M, \nabla(\lambda)) \not\cong 0$ for any $i$,
then $\lambda \in \Pi$.
This can be seen by noting that if 
$\Ext_S^i(M, \nabla(\lambda)) \not\cong 0$ then $M$ must contain a
composition factor $L(\mu)$ with $\mu \ge \lambda$. But since $M \in
\Mod (S(\Pi))$, $\mu \in \Pi$ and hence by saturation of $\Pi$,
$\lambda$ is in $\Pi$.

Thus if $M \in \Mod(S(\Pi))$ and using the isomorphism 
$\Ext^i_{S(\Pi)}(M,N) \cong \Ext^i_S(M,N)$,
$\wfd(M)$ as an $S(\Pi)$-module is
the same as $\wfd(M)$ as an $S$-module and similarly for the 
$\nabla$-filtration dimensions.
Hence if a module is unchanged by this form of truncation then its
$\Delta$- and $\nabla$-filtration dimensions are also unchanged.

Note however that if $M$ is an $S(\Pi)$-module then 
the injective and projective dimensions as a module for $S$ are
usually larger than those as a module for $S(\Pi)$.

\subsection{}\label{sect:trunc2}
Let $e=e_{\Gamma}$ as before, then the algebra $eSe$ is also quasi-hereditary,
with respect to $(\Gamma, \leq)$, and with standard modules
$e\Delta(\lambda)$ and costandard modules $e\nabla(\lambda)$, for
$\lambda \in \Gamma$ (see \cite[A3.11]{donkbk} or
\cite[\S1.6]{erdkluwer}). 
We note that $e \Delta(\mu) \ne 0 $ if and only if $\mu \in \Gamma$.
If $\mu$ is  in $\Gamma$
and $N$ is an $S$-module then
$\Ext_S^i(M,N) \cong \Ext^i_{eSe}(eM, eN)$ (\cite[proposition
  A3.13]{donkbk}).
%
%
Thus if $\mu \in \Gamma$ then the projective dimension of
$\Delta(\mu)$ 
is unchanged under this form of truncation.
But the $\Delta$- and $\nabla$-filtration dimensions usually are smaller for
$eSe$ than for $S$.

\subsection{}The two types of truncation are related by Ringel duality,
see \cite[theorem 3.4.6]{CPSstrat} or \cite[A4.9]{donkbk}.
We have that $eSe$ is a Ringel dual of a quasi-hereditary
quotient $S'/S'\epsilon S'$ where $S'$ is the Ringel dual
of $S$ and $\epsilon = \epsilon_{\Gamma}$ is defined as $e_\Gamma$ but
for $S^\prime$. Note that
$(\Gamma, \leq^{op})$ is a saturated subset
of the poset $(\partn, \leq^{op})$ for $S'$. 
So the conclusions in \ref{sect:trunc2}
also follow from \ref{sect:trunc1} together
with Proposition \ref{propn:ringel}.

\begin{example}
In the case of the Schur algebra $S(2,r)$ for $\GL_2$ we can now
completely describe the values of all the dimensions mentioned for
$L(\lambda)$, $\nabla(\lambda)$ and $\Delta(\lambda)$ under the
various forms of truncation and under Ringel duality.
For $S(2,r)$ the  poset $\partn(2,r)$ is totally ordered. Moreover
once we split the poset into block components then the resulting order
is a minimal one. 

Let $S$ be a block of $S(2,r)$ or of its Ringel dual $S^\prime(2,r)$ or
of any algebra obtained from $S(2,r)$ or of $S^\prime(2,r)$ by the two
forms of truncation defined above.

Suppose $S$ has $n+1$ simple modules. 
Then the simple modules for $S$ can be labelled by the numbers $0, 1, 2, 
\ldots, n$ with the usual ordering. We have
$\gfd(L(i))=\wfd(L(i))=\wfd(\nabla(i))=\gfd(\Delta(i))=i$ and
$\inj(\nabla(i))=\proj(\Delta(i))= n+1-i$.
We can also say that $\glob(S)=2\gfd(S)=2n$.
\end{example}

\section{Relating the partial order and inequalities \break for the
  homological dimensions}\label{sect:props}
\subsection{}
Let $S$ be a quasi-hereditary algebra with a duality fixing the simple modules.
We consider the following properties, that $S$ may not satisfy in
general but which are motivated 
by properties of the Schur algebras. 


We assume that the partial order $<$ is minimal in the sense of
section \ref{sect:min},
and that $S$ consists of one block (for $C$).

($A$) For all $\lambda, \mu \in \partn$,
if $\mu < \lambda$ then $\gfd(L(\mu))\leq \gfd(L(\lambda))$

($B$) For all $\lambda, \mu \in \partn$, 
if $\mu < \lambda$ then $\wfd(\nabla(\mu)) \leq \wfd(\nabla(\lambda))$

($C$) For all $\lambda \in \partn$, 
$\wfd(\nabla(\lambda)) = \gfd(S) - \inj(\nabla(\lambda))$.

($D$) For all $\lambda \in \partn$, $\wfd(\nabla(\lambda)) = \wfd(L(\lambda))$

($E$) For all $\lambda, \mu \in \partn$,
if $\mu < \lambda$ then $\inj(\nabla(\mu)) \geq \inj(\nabla(\lambda))$.

\subsection{Comparisons}  
First we list some easy observations. 
\begin{enumerate}
\item[$\bullet$]{We have $A$ implies $D$;
this follows by induction on $\leq$, using \cite[2.5]{parker1}.
Moreover, if 
$D$ holds then $A$ and $B$ are equivalent. So 
$A$ (and $D$) imply $B$.
}
\item[$\bullet$]{
$A$ and $C$ imply $E$.
}
\item[$\bullet$]{
If $C$ holds then $B$ and $E$ are equivalent.
}
\item[$\bullet$]{
By proposition \ref{propn:ringel} and its corollary \ref{cor:ringel}
we have:

$\bullet$ $B$ holds for $S$ if and only if $E$ holds for $S^\prime$
and

$\bullet$ $C$ holds for $S$ if and only if $C$ holds for $S^\prime$.
}

\end{enumerate}

\subsection{Examples} 

We now give a few examples which show that some of the reverse
implications do not hold.
\begin{example}
This shows that $A$, $B$ and $D$ do not imply $E$ or $C$.

We know that $A$, $B$, $C$, $D$ and $E$ all hold for the blocks of the
Schur algebra consisting of regular weights \cite{parker2}.
We also know that in the case $n=3$ (the first value of $n$ for which
there are primitive non-regular weights) that $A$, $B$ and $D$ hold
but that $C$ and $E$ do not hold in general for the non-regular
blocks.

As the representation theory for the Schur algebra is controlled by
that of the Special linear group
we now give an example for $\SL_3$  where
condition $E$ fails (and necessarily $C$ fails as well).
We use the standard notation and terminology of algebraic groups as in 
\cite{donkbk}.

Consider $S(3,6)$ for characteristic $2$. (A similar example works for
general characteristic.)  
The weights (in $\SL_3$-notation) of the non--simple block are
$$(0,0), \ (3,0), \ (0,3), \ (2,2), \ (4,1),\  (6,0)
$$
We observe that $(3,0) < (2,2)$, and we claim that
$\inj(\nabla(3,0))=1$ but $\inj(\nabla(2,2))=2$.

First, the injective $I(2,2)$ has a $\nabla$-filtration with quotients
$\nabla(2,2)$ and $\nabla(4,1)$ (only), this follows by reciprocity
from the decomposition matrix.  
Moreover, $I(4,1)$ has $\nabla$-quotients $\nabla(4,1), \nabla(6,0)$ and
$\nabla(6,0)$ is injective. This implies that $\nabla(2,2)$ has minimal
injective resolution
$$0 \to \nabla(2,2)\to I(2,2)\to I(4,1)\to I(6,0)\to 0.
$$
We claim now that $I(3,0)$ is isomorphic to the tilting module
$T(4,1)$ and hence is also projective. 

Now the tilting module $T(2,1)$ is isomorphic to $T(1,0) \otimes \St$
and hence is isomorphic to the injective hull of $L(1,0)$ as a
$G_1$-module. Thus it is indecomposable as a $G_1$ module so
by \cite[proposition 2.1]{donktilt}, 
the tilting module $T(4,1)$ is isomorphic to 
$T(2,1)\otimes T(1,0)^F$. The module $T(4,1)$ has simple socle
$L(3,0)$
and its injective hull is $I(3,0)$ and it follows that
$T(4,1) \cong I(3,0)$ since both have the same $\nabla$-quotients.

This implies that $I(3,0)/\nabla(3,0)$ is indecomposable. (It has
simple head $L(3,0)$.)
But ${\rm Ext}^1(\nabla(4,1),\break \nabla(2,2)) \cong k$ and 
we know that $I(2,2)$ is the non-split extension of $\nabla(4,1)$ and
$\nabla(2,2)$, so $I(3,0)/\nabla(3,0) \cong I(2,2)$.
Thus $\inj(\nabla(3,0))=1$.

Note that since $B$ holds but $E$ does not, this means that in
particular that $S(3,6)$ is not isomorphic to its Ringel dual.
Also note that the Ringel dual of $S(3,6)$ has property $E$ but not
$B$.
\end{example}

\begin{example}
This shows $B$ and $E$ do not imply $A$ or $C$ or $D$.

Let $S$ be the algebra $k{\QQ}/I$ where
${\QQ}$ is the quiver
$$\xymatrix@C=30pt{  
{\underset{\textstyle 0}\bullet} \ar@<2ex>@{->}[r]^{\textstyle\alpha_0}&
{\underset{\textstyle 1}\bullet} 
\ar@<2ex>@{->}[r]^{\textstyle\beta_0}\ar@<0ex>@{->}[l]^{\textstyle\alpha_1}&
{\underset{\textstyle 2}\bullet} 
\ar@<2ex>@{->}[r]^{\textstyle\gamma_0}\ar@<0ex>@{->}[l]^{\textstyle\beta_1}&
{\underset{\textstyle 3}\bullet} \ar@<0ex>@{->}[l]^{\textstyle\gamma_1}&
}
$$
with relations (composing on the right):
$$\alpha_1\alpha_0=0, \ \ \gamma_1\gamma_0=0, 
\ \ \beta_1\alpha_1=0=\beta_1\beta_0, \ \ \alpha_0\beta_0=0. 
$$

This is quasi-hereditary, with respect to the natural order on
$$\partn = \{ 0, 1, 2, 3 \}
$$
The composition factors of the $\nabla(i)$ are as follows
$$\begin{array}{r|cccc}
         &L(0)&L(1)&L(2)&L(3)\\
\hline
\nabla(0)& 1&  && \cr
\nabla(1)& 1&1& & \cr
\nabla(2)& 0&1&1& \cr
\nabla(3)& 0&1&1&1
\end{array}
$$
with  $\nabla(3)/L(3) \cong \nabla(2)$.
Then $2<3$ but $\gfd(L(3))=1$ and $\gfd(L(2))=2$. So this does not
satisfy $A$. Properties $D$ and $C$ also fail.
We have $\wfd(\nabla(3))=\wfd(\nabla(2))=2$, but $\inj(\nabla(2))=1$.
It does satisfy $B$ and $E$.
\end{example}

%
%
%
%

\begin{example} 
This shows that $D$ does not imply $A$ or $C$.

Let $S$ be the algebra $k{\QQ}/I$ where
${\QQ}$ is the quiver
$$\xymatrix@C=60pt@R=40pt{
{\underset{\textstyle 0}\bullet} \ar@<1.6ex>@{->}[r]^{\textstyle\alpha_0}
\ar@<2.4ex>@/^1pc/@{->}[r]^{\textstyle\epsilon_0} 
\ar@<0.8ex>@{->}[d]^{\textstyle\delta_1}&
{\underset{\textstyle 1}\bullet}
\ar@<0.8ex>@{->}[d]^{\textstyle\beta_0}\ar@<0ex>@{->}[l]^{\textstyle\alpha_1}
\ar@<0.6ex>@/^1pc/@{->}[l]^{\textstyle\epsilon_1}
\\
{\underset{\textstyle 3}\bullet} \ar@<1.6ex>@{->}[r]^{\textstyle\gamma_1}
\ar@<0.8ex>@{->}[u]^{\textstyle\delta_0}&
{\underset{\textstyle 2}\bullet}
\ar@<0ex>@{->}[l]^{\textstyle\gamma_0}\ar@<0.8ex>@{->}[u]^{\textstyle\beta_1}
}
$$
with the following relations:
\begin{align*}
\alpha_1\delta_1&=0=\alpha_1\epsilon_0 = \alpha_1\alpha_0,&
\beta_1\epsilon_1 &= \gamma_0\delta_0,&  
\gamma_1\gamma_0&=0=\beta_1\beta_0,
\\ 
\delta_0\delta_1&=0=\delta_0\epsilon_0=\delta_0\alpha_0,&
\epsilon_0\beta_0&=\delta_1\gamma_1,&
\epsilon_1\delta_1&=0=\epsilon_1\alpha_0 = \epsilon_1\epsilon_0. 
\end{align*}

%
%

We take the natural order on the index set $\partn = \{ 0, 1, 2, 3 \}$. 
The structure of the standard modules is as follows.

$$\rad(\Delta(1)) = \Delta(0)^2, \ \ \rad(\Delta(2)) = {\mathcal U}(L(1), L(0))
$$
$$\rad(\Delta(3)) = \Delta(0)\oplus \Delta(2).
$$
(writing ${\mathcal U}(-,-)$ for a uniserial module, listing 
composition factors starting at the top). We have a projective cover
$$0 \to P(2) \to P(1) \to \Delta(1) \to 0
$$
%
%
%
%
%
%
%

We have $\wfd(L(i)) = i$ for $i\leq 2$ and $\wfd(L(3))=1$, hence
property $A$ fails. 
We will now show that $\gfd(\Delta(i)) = i$ for $i\leq 2$ 
and $\gfd(\Delta(3))=1$
i.e. property $D$ holds.

This is clear for $i=0, 1$. Consider $i=2$.
We have the exact sequence 
$$0 \to \mathcal{U}(1, 0)\to \Delta(2)\to L(2)\to 0.
$$
The kernel has $\nabla$-filtration dimension equal to one 
and the cokernel has $\nabla$-filtration dimension equal to two. It follows
from the dual version of \cite[lemma 2.5(i)]{parker1} 
that $\gfd(\Delta(2))=2$. 
But we will need the following more explicit information about the
$\Ext$ groups.

\begin{lem}We have
$$\Ext^1(\Delta(1), \Delta(2)) = k, \ \ \Ext^2(\Delta(1), \Delta(2))=0; 
\leqno{(1)}$$
$$ \Ext^1(\Delta(0), \Delta(2))=k^2, \ \ 
\Ext^2(\Delta(0), \Delta(2)) = k.
\leqno{(2)}
$$
\end{lem}
\begin{proof} (1) From the projective cover for $\Delta(1)$ (see above) we get
$\Ext^i(\Delta(1), \Delta(2))$.

\noindent (2) 
We have $\Ext^2(\Delta(0),\Delta(2))\cong \Ext^2(\Delta(0),L(2))$
using $\gfd(\rad(\Delta(2)))=1$. 
This latter $\Ext$ group is isomorphic to
$\Ext^1(\Delta(0),Q)$ where $Q$ is the quotient $\nabla(2)/L(2)$.
This $\Ext$ group is one-dimensional as any non-split extension must
have simple socle $L(1)$ and hence must embed in $\nabla(1)$. Thus any
non-split extension is isomorphic to $\nabla(1)$ (by dimensions).

Clearly $\Hom(\Delta(0),\Delta(2))=k$, 
so we can now use the fact that 
$$\sum_{i} (-1)^i \Ext^i(\Delta(0),\Delta(2)) = 0$$
(see \cite[p.71]{ringel})
and that $\gfd(\Delta(2))=2$ to deduce that $\Ext^1(\Delta(0),\Delta(2))=k^2$.
%
%
\end{proof}

Now we will show that $\gfd(\Delta(3))=1$. 
That is we show that $\Ext^2(\Delta(i), \Delta(3))=0$ for $i\leq 2$. 

Apply $\Hom(\Delta(i), -)$ to the exact sequence
$$0 \to \Delta(2) \oplus \Delta(0) \to \Delta(3)\stackrel{\beta}\to L(3)\to 0
$$
Consider the resulting long exact sequence, for $i\leq 2$. 
Noting that $\Ext^j(\Delta(i), \Delta(0))=0$ for $i \ge 1$,
$\Hom(\Delta(i),L(3))=0$
and that $\gfd(L(3)=1$ this gives
\begin{multline*}
0 
\to \Ext^1(\Delta(i), \Delta(2))
\to \Ext^1(\Delta(i), \Delta(3))
\to \Ext^1(\Delta(i), L(3))
\\
\to \Ext^2(\Delta(i), \Delta(2))
\to \Ext^2(\Delta(i), \Delta(3))
\to 0.
\end{multline*}
Thus for $i=1$ or $2$ we have 
$\Ext^2(\Delta(i), \Delta(3))=0$. 

It remains to consider $i=0$, we substitute the dimensions proved in the
Lemma and get an exact sequence
$$0 \to k^2 \to \Ext^1(\Delta(0), \Delta(3))\stackrel{\beta^*}\to
k\to k \to \Ext^2(\Delta(0), \Delta(3))\to 0
$$
To complete the proof we will show that the map $\beta^*$ is zero. 

Take an element in $\Ext^1(\Delta(0), \Delta(3))$, say $\eta$, which is
represented by
$$0 \to \Delta(3)\to V\stackrel{\pi}\to \Delta(0)\to 0
$$
then $\beta^*(\eta)$ is represented by the push-out of $\beta$. Suppose
this is non-zero, then the middle term of the sequence
$\beta^*(\eta)$ must be uniserial with a simple top $L(0)$ and 
therefore, the top of $V$ must also be simple, isomorphic to $L(0)$.
So there is an epimorphism $\psi:\ e_0S = P(0) \to V$. 
Now we will use the relation $\delta_1\gamma_1 = \epsilon_0\beta_0$
to derive a contradiction.

We have $\psi(\epsilon_0) = 0$
(since $L(1)$ does not occur in the top of $\rad(V)$, i.e. of $\Delta(3)$). 
It follows that $\psi(\epsilon_0\delta_0) = \psi(\epsilon_0)\delta_0 = 0$.
Therefore also $\psi(\delta_1\gamma_1) = 0$. But on the other hand,
$\psi(\delta_1)$ generates $\ker(\pi) = \Delta(3)$ and
$\psi(\delta_1\gamma_1) = \psi(\delta_1)\gamma_1$ which spans 
$\Delta(3)e_2$ and is therefore non-zero,
a contradiction.
\end{example}

\subsection{}
We expect that $C$ is independent of any of the other conditions.
It certainly is not implied by any of them.
But to construct an example with $C$ but not $A$ 
seems to need rather a lot of technical detail
and would require many simple modules.

In summary $B$, $D$ and $E$ are independent of each other and none of
these imply $A$ or $C$, while $A$ implies $B$ and $D$ but not $C$ or $E$.

\subsection{} We now consider an algebra satisfying property $A$ and show
that the first condition of lemma \ref{lem:gfd} holds for this algebra.

\begin{example} \label{eg:propA}
In this example we consider a quasi-hereditary algebra $S$ with
duality preserving the simples that satisfies property $A$. We show
that $\glob(S)= 2\gfd(S)$ for this algebra.

Let $S$ be a quasi-hereditary algebra with four 
simple modules, and assume the $\nabla$-filtration dimensions of the
simples
$L(0)$, $L(1)$, $L(2)$, $L(3)$
are respectively $0, 1, 1, 2$. Then the quasi-hereditary quotient 
$S_2$ of $S$ obtained by factoring out $Se_3S$, where $e_3S = P(3)$, 
 must have 
$\gfd(S_2)=1 \le \glob(S_2) \le 2=2\gfd(S_2)$.

We would like to show that $\Ext_{S}^{4}(L(3), L(3)) \neq 0$.
To do this we need $\Ext_{S}^2(Q^\circ,Q)\cong \Ext_{S_2}^2(Q^\circ,
Q) 
\break\ne 0$ 
where $Q$ is the quotient $\nabla(3)/L(3)$.

We now show that if 
$Q$ is any $S_2$-module with $\gfd(Q)=1$ 
then $\Ext_{S_2}^2(Q^\circ, Q) \neq 0$. 


{\sc Case 1} $L(2)$ is not a composition factor of $Q$. Let $S_1$ be the
quasi-hereditary quotient of $S_2$ obtained by factoring $S_2e_2S_2$. (So that
$S_1 \in \mathcal{S}_1$ where $\mathcal{S}_1$ is as in section \ref{sect:glob}.) 
Then 
$Q$ is an $S_1$-module with $\gfd(Q)=1$ and then by theorem
\ref{thm:glob}  
we know that $\Ext_{S_1}^2(Q^\circ, Q)\cong
\Ext^2_{S_2}(Q^\circ,Q) \neq 0$.

{\sc Case 2 } The socle of $Q$ has only $L(2)$ as a composition factor.
Then consider the injective hull
$$0 \to Q \to I \to R \to 0
$$
where $I$ is isomorphic to a direct sum of copies of $\nabla(2)$.
Note that $R\neq 0$.
We have $\gfd(R)=0$ (by dimension shift).
That is, $R$ has a $\nabla$-filtration. Moreover $R$ is
an $S_1$-module. By dimension shift we get
$$\Ext_{S_2}^2(Q^\circ, Q) \cong \Ext_{S_2}^1(Q^\circ, R).$$

Now apply $\Hom_{S_2}(-, R)$ to the exact sequence
$$0 \to R^\circ\to I^\circ\to Q^\circ\to 0.
$$
This gives the exact sequence
$$0 \to \Hom_{S_2}(Q^\circ, R) 
\to \Hom_{S_2}(I^\circ, R)\to \Hom_{S_2}(R^\circ, R) 
\to \Ext_{S_2}^1(Q^\circ, R) \to 0.
$$
Since $L(2)$ is not a composition factor of $R$
and all the top composition factors of $I^\circ$ are $L(2)$
we deduce that $\Hom_{S_2}(I^\circ, R)=0$ and hence
the required $\Ext$-space is isomorphic to $\Hom_{S_2}(R^\circ, R)$.
This is certainly non-zero. 

{\sc General Case } Let $U$ be the largest submodule of $Q$ which does
not have $L(2)$ as a composition factor, and let $V$ be the quotient.
Then the socle of $V$ has only $L(2)$ as composition factors.
Consider the exact sequence of $S_2$-modules
$$0 \to U \to Q \to V \to 0.
$$
So each term has $\nabla$-filtration
dimension at most one.
Now $\gfd(Q)=1$, therefore at least one of $U$ and $V$ has 
$\nabla$-filtration dimension equal to one.

If $\gfd(V)=1$ then proceed as in the previous proof, to show that
there is a surjection 
$$\Ext_{S_2}^{2}(Q^\circ, Q) \rarr \Ext_{S_2}^2(V^\circ, V)$$
the latter $\Ext$ group being non-zero by Case 2. 

So assume now that $\gfd(V)=0$. Then actually $V$ is a direct sum of
$\nabla(2)$'s and hence is injective. Using this we show
that 
$$\Ext_{S_2}^2(Q^\circ, Q) \cong \Ext_{S_2}^2(U^\circ, U).
$$
But $U$ is as in Case 1 and hence this is non-zero.
\end{example}

Thus if property $A$ holds for a quasi-hereditary algebra with simple
preserving duality then if the global dimension is not twice the 
$\nabla$-filtration dimension then we must have at least five simples.

%
%

\providecommand{\bysame}{\leavevmode\hbox to3em{\hrulefill}\thinspace}
\providecommand{\MR}{\relax\ifhmode\unskip\space\fi MR }
\providecommand{\MRhref}[2]{%
  \href{http://www.ams.org/mathscinet-getitem?mr=#1}{#2}
}
\providecommand{\href}[2]{#2}

\end{document}